\documentclass[12pt]{amsart} 
\usepackage{amscd,amssymb,amsmath,amsthm,epsfig}
\usepackage[all]{xy}

\makeatletter
 \def\rightharpoonupfill@{\arrowfill@\relbar\relbar\rightharpoonup}
 \newcommand{\vect}{%
 \mathpalette{\overarrow@\rightharpoonupfill@}}
\makeatother

\newtheorem{theorem}{Theorem}[section]

\newtheorem{lemma}[theorem]{Lemma}

\theoremstyle{definition}
\newtheorem{remark}[theorem]{Remark}
\newtheorem{definition}[theorem]{Definition}
\newtheorem{example}[theorem]{Example}

\begin{document}
\title{Minimal Betti Numbers}

\thanks{This material is partially based upon work supported by the
  National Science Foundation under Grant No. DMS-0353622}

\author[C. Dodd]{Christopher Dodd}
\address{Christopher Dodd\\
2814 Rittenhouse St., NW\\
Washington, D.C. 20015}

\email{csdodd@gmail.com} 

\author[A. Marks]{Andrew Marks}
\address{Andrew Marks\\ 540 Linda Falls Terrace\\ Angwin, CA 94508}
\email{asmarks@gmail.com}

\author[V. Meyerson]{Victor Meyerson}
\address{Victor Meyerson\\ 19520 Cohasset Street\\ Reseda, CA 91335}
\email{vmeyerso@calpoly.edu}

\author[B. Richert]{Ben Richert$^1$}\thanks{$^1$Corresponding author:
  brichert@calpoly.edu} \address{Ben Richert\\ Mathematics
  Department\\ Cal Poly\\ San Luis Obispo, CA 93407}
\email{brichert@calpoly.edu}

\maketitle

\begin{abstract}
  We give conditions for determining the extremal behavior for the
  (graded) Betti numbers of squarefree monomial ideals. For the case
  of non-unique minima, we give several conditions which we use to
  produce infinite families, exponentially growing with dimension, of
  Hilbert functions which have no smallest (graded) Betti numbers
  among squarefree monomial ideals and all ideals. For the case of
  unique minima, we give two families of Hilbert functions, one with
  exponential and one with linear growth as dimension grows, that have
  unique minimal Betti numbers among squarefree monomial ideals.
\end{abstract}
\section{Introduction}

Let $R$ be a polynomial ring over a field $K$. Then given a Hilbert
function $H$ it is easy to see that there can be more than one ideal
$I\subset R$ such that the Hilbert function of $R/I$ is $H$. One can
further distinguish such ideals by passing to a finer invariant, the
graded Betti numbers, which gives rise to the question: given a
particular Hilbert function $H$, what sets of graded Betti numbers
actually occur? That this problem is bounded above, and hence finite,
is due to an important result by Bigatti and Hulett
\cite{Bigatti,Hulett} (independently in characteristic zero), and
Pardue \cite{Pardue} (in characteristic $p$) which says that, given a
Hilbert function $H$, the lexicographic ideal attaining $H$ has
everywhere largest graded Betti numbers. In fact, this says that the
partial order on the set of sets of graded Betti numbers of ideals
attaining a given Hilbert function has a unique maximum element.
Shortly thereafter, it was shown by Charalambous and Evans
\cite{CharalambousEvans} that this order need not have a unique
minimal element. Examples of infinite families of Hilbert functions
which did not support unique minimal elements were given by Richert
\cite{Richert}.

The structure of the partial order on graded Betti numbers has also
been considered on interesting subsets of the set of all ideals
attaining a given Hilbert function. For instance, Geramita, Harima,
and Shin \cite{GeramitaHarimaShin} showed that if one restricts to the
graded Betti numbers arising from a certain (dense) set of Gorenstein
ideals, then the associated partial order has a unique maximal element
which is constructed using a lexicographic ideal (Migliore and Nagel
\cite{MiglioreNagel} were later able to extend this result to an even
larger set of Gorenstein ideals attaining a given Hilbert function)
while Richert \cite{Richert} showed that a unique minimal element need
not exist. For stable ideals in dimension at most three, Francisco
\cite{Francisco} showed that, unlike the general and Gorenstein cases,
there is always a unique smallest element.

It is a result of Aramova, Herzog, and Hibi \cite{AramovaHerzogHibi1,
  AramovaHerzogHibi2} that if one restricts to the graded Betti
numbers arising from squarefree monomial ideals attaining a given
Hilbert function, then the associated partial order has a unique
maximal element (which arises, not surprisingly, from the squarefree
lexicographic ideal). Squarefree monomial ideals are particularly
interesting (to algebraists, topologists, and combinatorists alike)
because to each squarefree monomial ideal $I$ in $n$ variables can be
associated a simplicial complex $\Delta_I$ on $n$ vertices, while the
Hilbert function of $R/I$ is related to the face counts of $\Delta_I$
and the graded Betti numbers of $R/I$ can be computed by considering
sums of the ranks of the reduced homologies of subsets of $\Delta_I$.
It was known that (Gelvin, LaVictore, Reed, Richert
\cite{GelvinLaVictoreReedRichert}) for $n\le 5$ variables (and after
fixing a finite field), the partial orders arising from fixing a
Hilbert function were totally ordered, but that this failed in six
variables where, in fact, there is an example of a partial order which
does not have a unique smallest element. In the current paper, we
continue this line of inquiry. First, we generate an infinite family
(the size of which grows exponentially with dimension) of Hilbert
functions for which the partial order on the graded Betti numbers
corresponding to squarefree monomial ideals fails to have a unique
minimal element. We are able to show that this same family of Hilbert
functions gives rise to partially ordered sets corresponding to all
ideals (not only squarefree monomial ideals) which fail to have a
unique minimal element. We then find an infinite family (again,
growing exponentially) of Hilbert functions for which the partial
order on the Betti numbers (not graded Betti numbers) of squarefree
monomial ideals fails to have a unique smallest element (and are again
able to show that this family gives posets without unique minimal
elements in the case of all ideals). Next, we find an infinite family
(growing exponentially) of Hilbert functions for which the partial
order associated to the graded Betti numbers of squarefree monomial
ideals has a single element, and thus a unique smallest element. We
note an analogous family in the general case.  Finally, we find a
infinite family (growing linearly) of Hilbert functions for which the
partial order associated to the graded Betti numbers of squarefree
monomial ideals has a unique smallest element, and a nontrivial poset
tree.

\section{Background}

Let $R = K[x_1, \ldots, x_n]$ where $K$ is a field.  In what follows,
all ideals will be homogeneous.  In this paper, we will be studying
the graded Betti numbers of ideals with a fixed Hilbert function.
Recall that given a homogeneous ideal $I \subset R$, the \emph{Hilbert
  function} of $R/I$ in degree $d$, denoted as $H(R/I,d)$ is given by
$$H(R/I,t) = \dim_K(R/I)_t.$$
Furthermore, we will write
$$
0 \longrightarrow \bigoplus_{j \in \mathbb{Z}}R(-j)^{\beta^I_{h,j}}
\longrightarrow \cdots \longrightarrow \bigoplus_{j \in
\mathbb{Z}}R(-j)^{\beta^I_{1,j}} \longrightarrow R \longrightarrow
R/I \longrightarrow 0
$$
to be a minimal free resolution of $R/I$, and
$\beta^I=\{\beta^I_{i,j}\}$ to be the set of graded Betti numbers of
$R/I$.  It is useful to display graded Betti numbers in the following
table known as a Betti Diagram (using the notation of Macaulay 2
\cite{M2}).
$$
\begin{array}{c|cccccc}
    & s_0 & s_1 & s_2 & s_3 & \ldots \\ \hline
  0 & \beta_{0,0} & \beta_{1,1} & \beta_{2,2} & \beta_{3,3} & \dots \\
  1 & \beta_{0,1} & \beta_{1,2} & \beta_{2,3} & \beta_{3,4} & \dots \\
  2 & \beta_{0,2} & \beta_{1,3} & \beta_{2,4} & \beta_{3,5} & \dots \\
  3 & \beta_{0,3} & \beta_{1,4} & \beta_{2,5} & \beta_{3,6} & \dots \\
  \vdots & \vdots & \vdots & \vdots & \vdots & \ddots \\
\end{array}
$$
where $s_i$ is the sum of the entries in the $i^{\text{th}}$
column.  Then $S = \{s_i\}$ is the set of Betti numbers. There is an
obvious partial order on the Betti diagrams which arise for ideals
with a given Hilbert function. If $\beta^I$ and $\beta^J$ are the
graded Betti numbers of the ideals $I$ and $J$, then we say that
$\beta^I \ge \beta^J$ if $\beta^I_{i,j} \ge \beta^J_{i,j}$ for all $i$
and $j$.  Furthermore, $\beta^I > \beta^J$ if $\beta^I \ge \beta^J$
and there is a pair $(i,j)$ such that $\beta^I_{i,j} > \beta^J_{i,j}$.
We are interested in the extremal properties of this ordering. A
useful definition in this regard is that of $q$-linearity: A minimal
free resolution of an ideal $I$ is called $q$-linear if $I$ is
minimally generated by $q$-forms and $\beta_{i,j}= 0$ for each $j\neq
q+i-1$ and $j\neq 0$.

It turns out that the graded Betti numbers are a finer invariant
than Hilbert functions.  They are related by the following useful
equation.
\begin{theorem}
\cite{Stanley} Given an ideal $I \subset R = K[x_1, \ldots, x_n]$
with graded Betti numbers $\beta^I$,
$$
\sum_{d=0}^\infty H(R/I,d)t^d = \frac{\sum_{j=0}^\infty \sum_{i=0}^n
(-1)^i \beta^I_{i,j} t^j}{(1-t)^n}
$$
\end{theorem}

It follows from this formula that the diagonal alternating sums of a
Betti diagram are invariant for all Betti diagrams of ideals
attaining a given Hilbert function. We define this alternating sum
as follows:
\begin{definition}
Given an ideal $I$, the $j^{\text{th}}$ diagonal alternating sum of
the Betti diagram of $I$ is: \[ d_{j}={\displaystyle
\sum_{i=0}^{j}}(-1)^{i}\beta_{i,j}\] where the $\beta_{i,j}$ are the
graded Betti numbers of $I$.
\end{definition}

The partial ordering on Betti diagrams has a unique maximal Betti
diagram \cite{Bigatti, Hulett, Pardue}. In order to identify an ideal
attaining the maximal graded Betti numbers, we need the following
notation: Let $x_1^{a_1}x_2^{a_2}\cdots x_n^{a_n}$ and
$x_1^{b_1}x_2^{b_2}\cdots x_n^{b_n}$ be monomials such that $\sum a_i
= \sum b_i$ (i.e. they are of the same degree). We say
$x_1^{a_1}x_2^{a_2}\cdots x_n^{a_n} >_{lex} x_1^{b_1}x_2^{b_2}\cdots
x_n^{b_n}$ if and only if the first nonzero entry of $(a_1 - b_1, a_2
- b_2, \ldots, a_n - b_n)$ is positive.  This is known as the
\emph{lexicographic} or \emph{lex} ordering on monomials. An ideal $L$
is a \emph{lex ideal} if for all monomials $m \in L_d$ and $n \in
R_d$, then $n \ge_{lex} m$ implies $n \in L_d$. A lex ideal achieves
the maximal graded Betti numbers among all ideals attaining its
Hilbert function.

In this paper we are mainly concerned with squarefree monomial
ideals. A \emph{squarefree monomial ideal} is an ideal in $R =
K[x_1,\ldots,x_n]$ minimally generated by elements of the form
$x_1^{a_1}x_2^{a_2} \cdots x_n^{a_n}$, where $a_i \in \{0,1\}$.

We will make use of the fact that the squarefree monomial ideals are
in one-to-one correspondence with simplicial complexes.  Recall the
definition of a simplicial complex: A \emph{simplicial complex}
$\Delta$ on the vertex set $\{1,\ldots,n\}$ is a collection of
subsets called faces or simplices, closed under taking subsets; that
is if $\sigma \in \Delta$ is a face and $\tau \subseteq \sigma$,
then $\tau \in \Delta$.  A simplex $\sigma \in \Delta$ of
cardinality $|\sigma| = i+1$ has dimension $i$ and is called an
$i$-face.

The (well known, see for instance \cite{Stanley}) procedure to pass
from a simplicial complex $\Delta$ to a squarefree monomial ideal is
to form the ideal generated by monomials corresponding to the minimal
non-faces of $\Delta$, after renaming the vertex set of $\Delta$ to be
$\{x_1,\dots,x_n\}$.
\begin{example}
  The simplicial complex
  $\Delta=\big\{\{1,4\},\{2,3\},\{2,4\},\{3,4\}\big\}$ corresponds
  to the ideal $I=(x_1x_2, x_1x_3, x_2x_3x_4)$. We use dashed lines in
  the picture to indicate minimal non-faces.
\begin{center}
\includegraphics{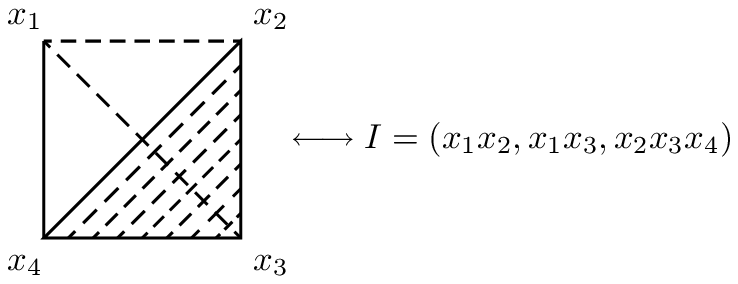}
\end{center}
\end{example}
This correspondence gives a bijection between simplicial complexes and
squarefree monomial ideals with no linear terms. In general, we will
conflate simplicial complexes and ideals due to this correspondence,
and so we are free to talk about the Betti numbers of a simplicial
complex or the homology of an ideal.

Given a simplicial complex with $f_i$ $i$-faces, its corresponding
\emph{$f$-vector} is the $n$-tuple $(f_0, f_1, \ldots, f_{n-1})$.  For
example, consider the simplicial complex \begin{center}
  \includegraphics{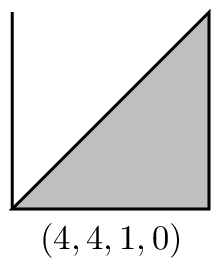}
\end{center} which has $f$-vector $(4,4,1,0)$. Here the $f$-vector
counts the 4 0-faces (or points), 4 1-faces (or edges), 1 2-face, and
0 3-faces (or volume).

The $f$-vector is an important invariant of a simplicial complex
because of the following:

\begin{theorem}
  \cite{Stanley} Let $I$ be the Stanley-Reisner ideal of a simplicial
  complex with $f$-vector $(f_0, f_1, \ldots, f_{n-1})$. Then
$$
H(R/I,m) =
  \begin{cases}
  1 & \text{$m = 0$,} \\
  \sum_{i=0}^{n-1}f_i \binom{m-1}{i} & \text{$m > 0.$}
  \end{cases}
$$
\end{theorem}

From the theorem, it follows that $f$-vectors and Hilbert functions
are in one-to-one correspondence.

Given a Hilbert function which arises for squarefree monomial ideals,
we will be particularly interested in the corresponding squarefree lex
ideal (because it always exists and is known to exhibit the maximal
Betti diagram \cite{AramovaHerzogHibi1, AramovaHerzogHibi2}). An ideal
$L$ is a \emph{squarefree lex ideal} if for squarefree monomials $m
\in L_d$ and $n \in R_d$, $n \ge_{\text{lex}} m$ implies $n \in L_d$.

A useful tool in studying simplicial complexes (and therefore their
associated ideals) is simplicial homology. If $\Delta$ is a simplicial
complex, we let $\Delta^{l}$ be the set of all $l$-faces in $\Delta$,
and let $C_{l}(\Delta)$ denote the $K$ vector space whose basis
consists of $\Delta^{l}$. We define the boundary map
$\partial_{l}:$$C_{l}(\Delta)\rightarrow C_{l-1}(\Delta)$to be \[
\partial_{l}(\{ i_{1},\ldots,i_{l+1}\})=\sum_{j=1}^{l+1}(-1)^{j+1}\{
i_{1},\ldots,\hat{i_{j}},\ldots,i_{l+1}\}\] where the hat indicates
omission. We define the $\Delta^{l}$ subspaces
$Z_{l}(\Delta)=\ker($$\partial_{l})$ and
$B_{l}(\Delta)=\textrm{image}(\partial_{l+1})$. By standard algebraic
topology, we have $B_{l}(\Delta)\subseteq Z_{l}(\Delta)$, and so we
can define the $l-homology$ of $\Delta$, $H_{l}(\Delta)$, to be the
$K$ vector space $Z_{l}(\Delta)/B_{l}(\Delta)$. In addition, the
$reduced$ $l-homology$ $\tilde{H_{l}}(\Delta)$ is given by
$H_{0}(\Delta)\cong\tilde{H_{0}}(\Delta)\oplus K$, and $H_{l}(\Delta)$
when $l>0$ . We now present the connection between the simplicial
homology of certain subsets of $\Delta$ and graded Betti numbers, in
the form of Hochster's formula. To properly express this formula, we
need a new notation. If $W\subseteq\{1,\ldots n\}$, then given a
simplicial complex $\Delta$, let $\Delta_{W}$ be the simplicial
complex defined by $\Delta_{W}=\Delta\cap P(W)$ where $P$ denotes
power set. Then we have

\begin{theorem}
  \cite{Hochster} Let $I\subseteq R$ be a squarefree monomial ideal,
  with $\Delta$ the associated simplicial complex. Then we have\[
  \beta_{i,j}^{I}=\sum_{W\subseteq\{1,\ldots,n\},|W|=j}
  \dim_{K}\tilde{H}_{j-i-1}(\Delta_{W};K)\] for all i and j.
\end{theorem}

\section{Non-unique Minimal Graded and Nongraded Betti Numbers}

\subsection{Graded Betti Numbers}

In this section we construct a fast-growing family of Hilbert
functions which fail to have unique minimal graded Betti numbers
among both squarefree monomial ideals, and all ideals. We proceed
using simplicial complexes, giving a method for preserving parts of
the Betti diagrams which guarantee incomparability. Our methods
allow us to double the number of $f$-vectors in each dimension with
this property, and so our family of $f$-vectors grows exponentially
with the dimension of the polynomial ring.

\begin{definition}
  Let $\Delta$ be a simplicial complex on $n$ vertices. We define the
  $j$-cone of $\Delta$ on the vertex $\{ n+1\}$, denoted
  $C_{(j)}\Delta$, to be the simplicial complex on $n+1$ vertices such
  that $\{ i_{1},\dots,i_{k}\}\in C_{(j)}\Delta$ if and only if
  either $\{ i_{1},\dots,i_{k}\}\in\Delta$ or $i_{k}=n+1$, $\{
  i_{1},\dots,i_{k-1}\}\in\Delta$ and $k-1\leq j$. As suggested by the
  definition, we define $C_{(\infty)}\Delta$ to be the $n$-cone of
  $\Delta$, $C\Delta$.
\end{definition}

\begin{example}
  Let $\Delta$ be the simplicial complex
  $$\Delta=\{\{1\}, \{2\}, \{3\},\{4\},\{1,2\}, \{1,4\}, \{2,4\},
  \{3,4\}, \{1,2,4\}\}\}\subset \{1,2,3,4\},$$ with $f$-vector
  $(4,4,1,0)$. 
  Then
\begin{eqnarray*}
C_{(0)}(\Delta)&=&\{\{1\}, \{2\}, \{3\},\{4\}, \{5\}, \{1,2\},
  \{1,4\}, \{2,4\}, \{3,4\}, \{1,2,4\}\}, \\
C_{(1)}(\Delta)&=&\{\{1\}, \{2\}, \{3\},\{4\}, \{5\}, \{1,2\},
  \{1,4\}, \{1,5\}, \{2,4\}, \{2,5\}, \{3,4\},\\&&\{3,5\}, \{4,5\},
  \{1,2,4\}\},\\
C_{(2)}(\Delta)&=&\{\{1\}, \{2\}, \{3\},\{4\}, \{5\}, \{1,2\},
  \{1,4\}, \{1,5\}, \{2,4\}, \{2,5\}, \{3,4\},\\ &&\{3,5\}, \{4,5\},
  \{1,2,4\}, \{1,2,5\}, \{1,4,5\}, \{2,4,5\}, \{3,4,5\}\}, \\
C_{(3)}(\Delta)&=&\{\{1\}, \{2\},
  \{3\},\{4\}, \{5\}, \{1,2\}, \{1,4\}, \{1,5\}, \{2,4\}, \{2,5\},
  \{3,4\},\\ &&\{3,5\}, \{4,5\}, \{1,2,4\}, \{1,2,5\}, \{1,4,5\},
  \{2,4,5\}, \{3,4,5\},\\& & \{1,2,4,5\}\},
\end{eqnarray*}
and, of course, $C_{(3)}(\Delta)=C_{(\infty)}(\Delta)$. Here $C_{(0)}$
has $f$-vector $(5,4,1,0,0)$, $C_{(1)}$ has $f$-vector $(5,8,1,0,0)$,
$C_{(2)}$ has $f$-vector $(5,8,5,0,0)$, and $C_{(3)}$ has $f$-vector
$(5,8,5,1,0)$.
\end{example}

On of the nice things about coning a simplicial complex is that the
graded Betti numbers do not change. We show below that we can
similarly preserve certain of the graded Betti numbers of a simplicial
complex after $j$-coning.

\begin{lemma}
\label{lem:diagonalsPreserved}Given a simplicial complex $\Delta$,
the Betti diagram of the simplicial complex $C_{(j)}\Delta$ will be
identical to the Betti diagram of $\Delta$ for the first $j+1$
diagonals.
\end{lemma}
\begin{proof}
  We employ Hochster's formula. Let $k\leq j+1$. Then we have \[
  \beta_{i,k}^{C_{(j)}(\Delta)}=
  \sum_{W\subseteq\{1,\ldots,n+1\},|W|=k}
  \dim_{K}\tilde{H}_{k-i-1}(\Delta_{W};K)\]
  \[=\sum_{W\subseteq\{1,\ldots
    n\},|W|=k}\dim_{K}\tilde{H}_{k-i-1}(\Delta_{W};K)+\sum_{\{
    n+1\}\in W,|W|=k}\dim_{K}\tilde{H}_{k-i-1}(\Delta_{W};K)\]
  However, from the definition of $j$-coning and the complex
  $\Delta_{W}$, we have that for any $W$ with $\{ n+1\}\in W$ and at
  most $j$ vertices in $\{1,\ldots n\}$, $\Delta_{W}$ is a cone. By
  standard algebraic topology, a cone is contractible and so all of
  its reduced homology spaces are zero. Therefore, the second term in
  the above sum vanishes, and as the first term is
  $\beta_{i,k}^{\Delta}$ by Hochster's formula, the lemma is proved.
\end{proof}
\begin{definition}
\label{def:coningdefs}Given $j$, and an $f$-vector
$(f_{0},f_{1},\dots,f_{n-1})$ of a simplicial complex $\Delta$, define
$f_{m}^{[k]}$ to be the $m^{\text{th}}$ entry of the $f$-vector
$(C_{(j)})^{k}\Delta$. Also, for $m>j$, define\[
f_{m}^{(k)}=\begin{cases}
  f_{m} & \text{for }m=j+1\text{ or }k=0\text{,}\\
  f_{m-1}^{(k-1)}+f_{m}^{(k-1)} & \text{otherwise.}\end{cases}\]
\end{definition}
\begin{remark}
  If the $f$-vector of $\Delta$ is $(f_{0},f_{1},\dots,f_{n-1})$, then
  it is easy to show that the $f$-vector of $C_{(j)}\Delta$ is: \[
\begin{cases}
  (1+f_{0},f_{0}+f_{1},f_{1}+f_{2},\dots,f_{j-1}+f_{j},f_{j+1},
  \dots,f_{n-1},0) & j<n\\
  (1+f_{0},f_{0}+f_{1},f_{1}+f_{2},\dots,f_{n-2}+f_{n-1},f_{n-1}) &
  j\geq n\end{cases}\] or, in the notation of Definition
\ref{def:coningdefs}, \[
 \begin{cases}
(f^{[1]}_{0},f^{[1]}_{1},\dots,f^{[1]}_{j},f_{j+1},\dots,f_{n-1},0) & j<n\\
(f^{[1]}_{0},f^{[1]}_{1},\dots,f^{[1]}_{n-1},f_{n-1}) & j\geq
n\end{cases}\]
 We speak of this $f$-vector as the one generated by $j$-coning
the initial $f$-vector. Consider the family of $f$-vectors derived
by starting with an initial $f$-vector and at each stage, both
$j$-coning and $\infty$-coning it. The first iteration of such a
tree is shown below, assuming $j<n$:
\[ {\tiny \xymatrix{
 & (f_{0},\dots,f_{n-1})\ar[dl]_{j\text{-cone}}\ar[dr]^{\infty\text{-cone}} &           \\
(f_{0}^{[1]},\dots,f_{j}^{[1]},f_{j+1},\dots,f_{n-1})\ar@{.}[d] & & (f_{0}^{[1]},\dots,f_{j}^{[1]},f_{j}^{[0]}+f_{j+1}^{(1)},\dots,f_{n-1}^{(1)})\ar@{.}[d]\\
& &}}\] We will use this technique to generate families of
$f$-vectors which grow exponentially with dimension.
\end{remark}
\begin{definition}
Given a simplicial complex $\Delta$ and a $k$-tuple
$(m_{0},m_{1},\dots,m_{k})$ where each $m$ is either a nonnegative
integer or $\infty$, define $C_{(m_{0},m_{1},\dots,m_{k})}\Delta$ to
be the simplicial complex generated by $m_{0}$-coning, then
$m_{1}$-coning, and so on. Similarly, if $I$ is the associated ideal
of $\Delta$, define $C_{(m_{0},m_{1},\dots,m_{k})}I$ to be the ideal
associated to $C_{(m_{0},m_{1},\dots,m_{k})}\Delta$.
\end{definition}
\begin{lemma}
\label{lem:distinctTree}Given an $f$-vector
$(f_{0},f_{1},\dots,f_{j},\dots,f_{c},\dots)$ of a simplicial
complex $\Delta$ where $f_{c}$ is the last nonzero element of the
$f$-vector and $j\leq c$, the $f$-vectors of the simplicial
complexes in the tree generated by repeatedly $j$-coning and
$\infty$-coning $\Delta$ are distinct.
\end{lemma}
\begin{proof}
  We can index each $f$-vector in the tree with the $k$-tuple
  $(m_{0},m_{1},\dots,m_{k})$ where each $m$ is either $j$ or $\infty$
  and $C_{(m_{0},m_{1},\dots,m_{k})}\Delta$ is the corresponding
  simplicial complex. Let $\{ t_{i}\}$, $0\leq i<r$ be the $t$ such
  that $m_{t_{i}}=\infty$ so that $r$ is the number of times we
  $\infty$-cone. Then the $f$-vector of
  $C_{(m_{0},m_{1},\dots,m_{k})}\Delta$ is: \[
  (f_{0}^{[k]},f_{1}^{[k]},\dots,f_{j}^{[k]},{\displaystyle
    \sum_{i=0}^{r-1}}f_{j}^{[t_{i}]}+f_{j+1}^{(r)},{\displaystyle
    \sum_{i=0}^{r-2}}f_{j}^{[t_{i}]}+f_{j+2}^{(r)},\]
  \[\dots,f_{j}^{[t_{0}]}+f_{j+r}^{(r)},f_{j+r+1}^{(r)},\dots,f_{c+r}^{(r)},\dots)\text{.}\]
  From this $f$-vector we can read $r$ from the index of the last
  nonzero entry minus $c$. Further, as the $f_{j}^{p}$ are
  monotonically increasing and thus distinct, the elements at position
  $j+1$ through $j+r$ determine the $\{ t_{i}\}$, and as the
  $j^{\text{th}}$ element gives $k$, the $f$-vector uniquely
  determines $(m_{0},m_{1},\dots,m_{k})$.
\end{proof}
\begin{remark}
This Lemma ensures that $j$-coning and $\infty$-coning a suitable
$f$-vector gives a family of Hilbert functions which grows as
$2^{n}$ with dimension. It can be shown that this growth is not
necessarily achieved if we $j$-cone and $l$-cone for arbitrary $j$
and $l$.
\end{remark}
\begin{theorem}
\label{theorem:firstfamily} Given a set of squarefree monomial ideals
$I_{k}$ corresponding to an $f$-vector $\vect{v}$, if for some $j$
their Betti diagrams have $j^{\text{th}}$ diagonal such that the
alternating sum $d_{j}$ is nonzero and, for each $i$, ${\displaystyle
  \min_k \beta_{i,j}^{I_{k}}=0}$, then each $f$-vector in the tree
created by $(j-1)$-coning and $\infty$-coning $\vect{v}$ has an
associated Hilbert function such that the poset tree of graded Betti
numbers of all ideals attaining this Hilbert function will fail to
have a unique minimum. The poset tree of each Hilbert function also
fails to have a unique minimum if we restrict to the graded Betti
numbers of all squarefree monomial ideals.
\end{theorem}
\begin{proof}
We show that this condition on the $j^{\text{th}}$ diagonal
guarantees that a unique minimum cannot exist among the graded Betti
numbers. As this diagonal is preserved under both $(j-1)$-coning and
$\infty$-coning and the $C_{(m_{0},\dots,m_{k})}I_{k}$ are
squarefree monomial ideals, the theorem follows. As ${\displaystyle
\min}\beta_{i,j}^{I_{k}}=0$ for all i, any diagram that is less than
all $\beta^{I_{k}}$ in the partial order has all zeros on its
$j^{\text{th}}$ diagonal. However, as $d_{j}$ is nonzero, no Betti
diagram can have all zeros along the $j^{\text{th}}$ diagonal and
attain $d_{j}$. Thus, no ideal can be less than all $\beta^{I_{k}}$
and attain this Hilbert function, and as the conditions on the $I_k$
imply that at least two of the $I_k$ are incomparable, a unique
minimum cannot exist.
\end{proof}
\begin{remark}
  The number of (not necessarily distinct) $f$-vectors in the tree
  constructed in Theorem \ref{theorem:firstfamily} grows as $2^{c}$
  where $c$ is the number of times the initial $f$-vector has been
  coned. If the initial $f$-vector satisfies the hypothesis of Lemma
  \ref{lem:distinctTree}, then we are guaranteed that these
  $f$-vectors will be distinct, and so this family grows exponentially
  with the dimension of the polynomial ring. We can enlarge this
  family by $(j-1),\dots,(j+c-1)$-coning as well as $\infty$-coning
  where $c$ is the number of times we have coned already, however, our
  computational evidence suggests that this family still grows
  exponentially.
\end{remark}
\begin{example}\label{e:ideals}
  From an exhaustive computational search, several examples were found
  that satisfy the hypothesis of the above theorem. The programs
  written to find these examples are available from the authors on
  request. The first example occurring in lex order was with the
  $f$-vector $(6,8,4,0,0,0)$, and had previously been noted by Gelvin,
  LaVictore, Reed, Richert \cite{GelvinLaVictoreReedRichert}. Two
  ideals attaining this $f$-vector are: \[
  I=(x_{1}x_{2},x_{1}x_{3},x_{2}x_{3},x_{3}x_{4},x_{3}x_{5},
  x_{3}x_{6},x_{4}x_{5})\]
\[
J=(x_{1}x_{2},x_{1}x_{4},x_{2}x_{3},x_{2}x_{5},
x_{3}x_{4},x_{4}x_{5},x_{4}x_{6}, x_{1}x_{3}x_{5}x_{6})\] with graded
Betti numbers (calculated over ${\mathbb{Z}}_{101}$ with Macaulay 2):
\begin{center}$\beta^{I}=$\begin{tabular}{c|cccccc} & 1& 7& 13& 11& 5&
    1\tabularnewline \hline 0& 1& 0& 0& 0& 0& 0\tabularnewline 1& 0&
    7& 12& 10& 5& 1\tabularnewline 2& 0& 0& 1& 1& 0& 0\tabularnewline
\end{tabular} \quad$\beta^{J}=$\begin{tabular}{c|ccccc}
& 1& 8& 14& 9& 2\tabularnewline \hline 0& 1& 0& 0& 0&
0\tabularnewline 1& 0& 7& 12& 8& 2\tabularnewline 2& 0& 0& 0& 0&
0\tabularnewline 3& 0& 1& 2& 1& 0\tabularnewline
\end{tabular} \end{center}The $6^{\text{th}}$ diagonal of these Betti diagrams
satisfies the condition that their alternating sums are not zero, but
the minimum in each position on the diagonal is zero. Furthermore,
$(C_{(\infty)})^{3}(6,8,4,0,0,0)$ has a nonzero entry in the
$5^{\text{th}}$ position of its $f$-vector. Thus, the family of
$f$-vectors created by $5$-coning and $\infty$-coning
$(C_{(\infty)})^{3}(6,8,4,0,0,0)$ has associated Hilbert functions
with non-unique minimal graded Betti numbers for both the squarefree
and general case and grows exponentially with dimension.
\end{example}

\subsection{Betti Numbers}

In this section we provide a method for constructing families of
Hilbert functions without unique minimal Betti numbers among both
squarefree monomial ideals, and all ideals. Again, we use methods
for preserving the properties of the Betti diagrams which guarantee
incomparability. Our family of $f$-vectors grows exponentially with
the dimension of the polynomial ring.

\begin{lemma}
\label{lem:diagonalLemma}If an ideal $I$ has the Betti diagram:
\begin{center}$\beta^{I}=$\begin{tabular}{c|cccc}
&
 $s_{0}$&
 $s_{1}$&
 $s_{2}$&
 $\dots$\tabularnewline
\hline $0$&
 $\beta_{0,0}$&
 $\beta_{1,1}$&
 $\beta_{2,2}$&
 $\dots$\tabularnewline
$1$&
 $\beta_{0,1}$&
 $\beta_{1,2}$&
 $\beta_{2,3}$&
 $\dots$\tabularnewline
$2$&
 $\beta_{0,2}$&
 $\beta_{1,3}$&
 $\beta_{2,4}$&
 $\dots$\tabularnewline
$\vdots$&
 $\vdots$&
 $\vdots$&
 $\vdots$&
 $\ddots$\tabularnewline
\end{tabular} \end{center}such that ${\displaystyle \sum_{i=0}^{j}}\beta_{i,j}=\mid{\displaystyle d_{j}}\mid$
for all $j$, then if L is an ideal with the same Hilbert function,
${\displaystyle \sum_{i}}s_{i}^{L}\geq{\displaystyle \sum_{i}}s_{i}$
. 
%In particular, $L$ does not have smaller Betti numbers or smaller graded Betti numbers than $I$.
\end{lemma}
\begin{proof}
  Given an ideal $L$ with the same Hilbert function as $I$,
  ${\displaystyle \sum_{j}s_{j}^{L}}={\displaystyle
    \sum_{i,j}}\beta_{i,j}^{L}\geq{\displaystyle \sum_{j}\mid
    d_{j}|}$.  As ${\displaystyle \sum_{j}s_{j}}={\displaystyle
    \sum_{j}\mid d_{j}|}$ for $I$, if $L$ had Betti numbers
  $s_{i}^{L}$ such that ${\displaystyle
    \sum_{i}}s_{i}^{L}<{\displaystyle \sum_{i}}s_{i}$, then
  ${\displaystyle \sum_{i}}s_{i}^{L}<\sum_{j}\mid d_{j}|$, a
  contradiction.
\end{proof}
\begin{lemma}
\label{lem:addPoint} Let $\Delta$ be a simplicial complex, and let
$I$ be its associated ideal. The graded Betti numbers of
$C_{(0)}\Delta$ can be computed as:
\[ \beta_{i,j}^{C_{(0)}\Delta}=\begin{cases}
  \beta_{i,j}^{I} & \text{for }i=0\text{,}\\
  \beta_{i-1,j-1}^{I}+\beta_{i,j}^{I}+\binom{n}{i} & \text{for }i=j-1\text{,}\\
  \beta_{i-1,j-1}^{I}+\beta_{i,j}^{I} & \text{otherwise.}\end{cases}\]
where $n$ is the number of vertices of $\Delta$.
\end{lemma}
\begin{proof}
  Recall that $C_{(0)}\Delta$ is the simplicial complex obtained by
  adding a point to $\Delta$. By definition, we have that
  $\tilde{H}_{l}(C_{(0)}\Delta)=\tilde{H}_{l}(\Delta)$ if $l>0$ and
  $\tilde{H}_{0}(C_{(0)}\Delta)=\tilde{H}_{0}(\Delta)+1$. By
  Hochster's formula:
\begin{eqnarray*}
  \beta_{i,j}^{C_{(0)}\Delta}&=&\sum_{W\subseteq\{1,\ldots
    n\},|W|=j}\dim_{K}\tilde{H}_{j-i-1}(\Delta_{W};K) \\ & &+\sum_{\{
    n+1\}\in W,|W|=j}\dim_{K}\tilde{H}_{j-i-1}(\Delta_{W};K)\\
&=&\beta_{i,j}^{I}+\sum_{\{ n+1\}\in
W,|W|=j}\dim_{K}\tilde{H}_{(j-1)-i}(\Delta_{W};K)\end{eqnarray*}

However, subsets of size $j$ containing $\{ n+1\}$ are in bijective
correspondence with subsets of $\{1,\ldots,n\}$ of size $j-1$. If
$j>i+1$, then by the above formula for adding a point we have that
each of the vector spaces in the last term above are isomorphic to
the corresponding spaces $\tilde{H}_{(j-1)-i}(\Delta_{W'};K)$ where
$W'=W\setminus\{ n+1\}$, and so the last sum is just
$\beta_{i-1,j-1}^{I}$. Finally, if $j=i+1$, then in addition to
$\beta_{i-1,j-1}^{I}$, we add $1$ for each subset of
$\{1,\ldots,n\}$ of size $j-1=i$, so this gives the second line of
the formula.
\end{proof}
\begin{lemma}\label{l:betti}
  Suppose that two squarefree monomial ideals $I$ and $J$ have the
  same $f$-vector $\vect{v}$ and have Betti numbers $s^{I}$ and
  $s^{J}$ such that for some $k$, $s_{0}^{I}=s_{0}^{J}$,
  $s_{1}^{I}>s_{1}^{J}$, $s_{k}^{I}<s_{k}^{J}$, $s_{k+i}^{I}\leq
  s_{k+i}^{J}$ for all positive $i$, and $\beta_{i,j}^{I}=0$ when
  $j-i$ is even and $j>0$. Then each $f$-vector in the tree created by
  $0$-coning and $\infty$-coning $\vect{v}$ has a Hilbert function
  such that the poset tree of all Betti numbers of ideals with this
  Hilbert function will not have a unique minimum.  This will remain
  true if we restrict to the Betti numbers of squarefree monomial
  ideals.
\end{lemma}
\begin{proof}
We first observe that the conditions of the theorem still hold under
$0$-coning and $\infty$-coning. Graded Betti numbers are unchanged
under $\infty$-coning, and using Lemma \ref{lem:addPoint}, it is
easy to show that the conditions still hold under $0$-coning. Thus,
for all choices $m_{i}\in\{0,\infty\}$, the Betti numbers of
$C_{(m_{0},\dots,m_{k})}I$ and $C_{(m_{0},\dots,m_{k})}J$ will be
incomparable and $C_{(m_{0},\dots,m_{k})}I$ will always satisfy
${\displaystyle
\sum_{j}s_{j}^{C_{(m_{0},\dots,m_{k})}I}}={\displaystyle
\sum_{j}\mid d_{j}^{C_{(m_{0},\dots,m_{k})}I}|}$. As Lemma
\ref{lem:diagonalLemma} proves that no ideal can have smaller Betti
numbers than $C_{(m_{0},\dots,m_{k})}I$, unique minimal Betti
numbers cannot exist for the family of $f$-vectors created by
$0$-coning and $\infty$-coning $\vect{v}$.
\end{proof}
\begin{remark}
  The theorem remains true with a nearly identical proof if we replace
  the clause $\beta_{i,j}^{I}=0$ if $j-i$ is an even number and $j>0$
  with $\beta_{i,j}^{J}=0$ if $j-i$ is an even number and $j>0$.
\end{remark}
\begin{example}
  Recall from example \ref{e:ideals} that the ideals
 \[ I=(x_{1}x_{2},x_{1}x_{3},x_{2}x_{3},x_{3}x_{4},x_{3}x_{5},
 x_{3}x_{6},x_{4}x_{5}),\ {\rm and}\]
\[
J=(x_{1}x_{2},x_{1}x_{4},x_{2}x_{3},x_{2}x_{5},
x_{3}x_{4},x_{4}x_{5},x_{4}x_{6}, x_{1}x_{3}x_{5}x_{6}),\]
corresponded to the $f$-vector $(6,8,4,0,0,0)$ and have graded Betti
numbers (over ${\mathbb{Z}}_{101}$)
\begin{center}$\beta^{I}=$\begin{tabular}{c|cccccc} & 1& 7& 13& 11& 5&
    1\tabularnewline \hline 0& 1& 0& 0& 0& 0& 0\tabularnewline 1& 0&
    7& 12& 10& 5& 1\tabularnewline 2& 0& 0& 1& 1& 0& 0\tabularnewline
\end{tabular} \quad$\beta^{J}=$\begin{tabular}{c|ccccc}
& 1& 8& 14& 9& 2\tabularnewline \hline 0& 1& 0& 0& 0&
0\tabularnewline 1& 0& 7& 12& 8& 2\tabularnewline 2& 0& 0& 0& 0&
0\tabularnewline 3& 0& 1& 2& 1& 0\tabularnewline
\end{tabular} \end{center}
These Betti numbers satisfying the conditions laid out in Lemma
\ref{l:betti} (with the roles of $I$ and $J$ reversed). Thus, the
family generated by $0$-coning and $\infty$-coning $(6,8,4,0,0,0)$
consists of $f$-vectors corresponding to Hilbert functions whose poset
trees of Betti numbers fail to have unique minimums. This family grows
exponentially with dimension by Lemma \ref{lem:distinctTree}.
\end{example}
\begin{remark}
The infinite family given above is also an infinite family of
Hilbert functions with non-unique minimal \emph{graded} Betti
numbers, as incomparable Betti numbers implies incomparable graded
Betti numbers.
\end{remark}

\section{Unique Minimal Betti Numbers}

\subsection{Squarefree Monomial Ideals}

One of the simplest ways to construct Hilbert functions for squarefree
monomial ideals with unique minimal graded Betti numbers is to find
those for which the graded Betti numbers of the squarefree lex ideal
are minimal. As the squarefree lex ideal always gives maximal graded
Betti numbers (among squarefree monomial ideals), uniqueness follows
immediately.

To proceed, we observe that a minimal free resolution of a squarefree
lex ideal generated in a single degree $d$ is $d$-linear (this was
proved by Aramova, Hibi, and Herzog \cite{AramovaHerzogHibi1}). In
fact, it is true (see Herzog, Reiner, and Welker in
\cite{HerzogReinerWelker}) that squarefree lex ideals are
componentwise linear (an ideal $I$ is componentwise linear if $I_d$ is
$d$-linear for all $d$). In particular, this implies that if $L$ is a
squarefree lex ideal with no minimal generators in degree $t$, then
$\beta^L_{i,j}=0$ for $j=t+i-1$---or in words, the Betti diagram of a
squarefree lex ideal $L$ may contain a nonzero entry in row $i$ only
if $i=0$ or $L$ has a minimal generator in degree $i+1$.

We now give a family, growing exponentially with dimension, of
$f$-vectors for which the corresponding poset tree of graded Betti
numbers (for squarefree monomial ideals) has a unique (and hence a
unique minimal) element.

\begin{theorem}
  Suppose that $L$ is a squarefree lex ideal generated in a single
  degree. Then the poset tree of graded Betti numbers of squarefree
  monomial ideals with the same Hilbert function as $R/L$ consists of
  a unique element (which is thus uniquely minimal). The family of
  $f$-vectors which give rise to such $L$ grows exponentially with
  dimension.
\end{theorem}
\begin{proof}
  We know that a squarefree monomial lex ideal, $L$, generated in a
  single degree has a $d$-linear resolution. Thus, by Lemma
  \ref{lem:diagonalLemma}, no other ideal with the same $f$-vector can
  have smaller Betti numbers than $L$ and so it is the unique minimum.
  For a polynomial ring with $n$ variables, there are $\binom{n}{i}$
  squarefree monomials of degree $i$ and thus $2^{n}-n-1$ distinct
  squarefree monomial lex ideals generated in a single degree. This
  gives the exponential growth with dimension.
\end{proof}
\begin{remark}
The same techniques can be used to show that Hilbert functions
containing lex ideals generated in a single degree will have unique
minimal graded Betti numbers among all ideals.
\end{remark}

\subsection{Unique Mins via Hochster's formula}\label{s:hochster}

We now present a new method for finding unique minima among
squarefree ideals using Hochster's formula.

\begin{lemma}
\label{lem:tworow} 
Suppose that $L$ is a squarefree lex ideal with $\beta^L_{i,j}\neq 0$
iff $i,j=0$ or $j-i=1\text{ or }2$. Then if an ideal exists that
attains the same $f$-vector as $L$, and satisfies the hypothesis of
Lemma \ref{lem:diagonalLemma}, it will have the unique minimal graded
Betti numbers among all squarefree monomial ideals with this
$f$-vector.
\end{lemma}
\begin{proof}
  The conditions on $L$ imply that any ideal with the same $f$-vector
  will have at most two nonzero Betti numbers in each of its diagonal
  alternating sums, and these will have opposite signs in the
  summation. Thus, if an ideal $I$ satisfies the hypothesis of Lemma
  \ref{lem:diagonalLemma}, it will have at most one nonzero graded
  Betti number on each diagonal and so no other ideal can have a
  smaller graded Betti number than $I$ while preserving the $d_{j}$.
\end{proof}
\begin{theorem}
\label{theorem:hochstermin} In the polynomial ring of dimension n,
the Hilbert function corresponding to the $f$-vector
$(n,k,0,\ldots,0)$ for $0\leq k\leq n$ has an ideal which attains
minimal graded Betti numbers among all squarefree monomial ideals with
the same Hilbert function.
\end{theorem}
\begin{proof}
  If $k=0$ we have that the given simplicial complex is the only one
  with the corresponding $f$-vector, and hence is minimal. For $k>0$,
  we note that the squarefree lex ideal for this $f$-vector can only
  be generated in degree 2 and 3, as a generator in degree 4 would
  correspond to a minimal non-3-face, which would certainly imply that
  this simplicial complex has 2-faces; but the above $f$-vector has
  none. By Lemma \ref{lem:tworow}, it remains to show that there exist
  ideals attaining $(n,k,0,\ldots,0)$ for each $1\le k \le n$ which
  satisfy the hypothesis of Lemma \ref{lem:diagonalLemma}. We do this
  in the following two results.
\end{proof}

In this first lemma, we show that for each $1\leq k \leq n-1$, there
is an ideal attaining $(n,k,0,\ldots,0)$ which satisfies the
hypothesis of Lemma \ref{lem:diagonalLemma}.

\begin{lemma}
  Let $1\leq k \leq n-1$. Then the simplicial complex $$\Delta=
  \{\{1\},\ldots,\{ n\},\{1,2\},\{2,3\},\ldots,\{ k,k+1\}\}$$
  has a
  $2$-linear resolution.
\end{lemma}
\begin{proof}
  By Hochster's formula, a resolution will be 2-linear if there is no
  reduced $l$-homology for $l\geq 1$, for any sub-complex $\Delta_{W}$
  (because in that case $\beta_{i,j}=0$ for $i=j$ and $\beta_{i,j}$,
  which for $j>i+1$ depends only on the spaces $\tilde{H}_{j-i-1}$,
  are all zero).  But note that for any $W\subseteq\{1,\ldots n\}$ the
  simplicial complex $\Delta_{W}$ is homotopy equivalent to a finite
  set of points. Since a finite set of points never has reduced
  $l$-homology for $l\geq 1$, the lemma is proved.
\end{proof}

We now demonstrate that there is an ideal attaining $(n,n,0,\ldots ,
0)$ which satisfies the hypothesis of Lemma \ref{lem:diagonalLemma}.

\begin{lemma}
  The only graded Betti numbers of the simplicial complex
\[ \{\{1\},\ldots,\{ n\},\{1,2\},\{2,3\},\ldots,\{ n,1\}\}\]
which can be nonzero are $\beta_{0,0}$, $\beta_{i,i+1}$ for $1\leq
i\leq n-2$ and $\beta_{n-2,n}$.
\end{lemma}
\begin{proof}
  The case $\beta_{i,j}$ for $j<n$ is done above as if
  $W\subset\{1,\ldots,n\}$ (strict containment) then $\Delta_{W}$ is a
  simplicial complex in the form of the above lemma. If $j=n$, then
  $\beta_{i,n}$ depends only the space $\tilde{H}_{n-i-1}(\Delta,K)$.
  However, this simplicial complex is homotopy equivalent to a circle,
  so its only nonzero reduced homology space is $\tilde{H}_{1}$, which
  occurs in Hochster's formula when $i=n-2$.
\end{proof}
\begin{remark}
  This family of Hilbert functions with unique minimal graded Betti
  numbers grows linearly with dimension; from Theorem
  \ref{theorem:hochstermin}, we have $n+1$ Hilbert functions in each
  dimension.
\end{remark}


\begin{thebibliography}{10}

\bibitem[A-H-H1]{AramovaHerzogHibi1} Aramova, A., Herzog, J., Hibi,
T., Squarefree lexsegment ideals, \emph{Mathematische Zeitschrift},
\textbf{228} (1998), 353--378.

\bibitem[A-H-H2]{AramovaHerzogHibi2} Aramova, A., Herzog, J., Hibi,
T., Shifting operations and graded Betti numbers, \emph{Journal of
Algebraic Combinatorics}, \textbf{12} (2000), 207--222.

\bibitem[B]{Bigatti} A. M. Bigatti, Upper bounds for the Betti numbers
    of a given Hilbert function, \emph{Comm. in Alg.}, \textbf{21}
    (1993), no.~7, 2317--2334.

\bibitem[C-E]{CharalambousEvans} H.~Charalambous and E.~G. Evans, Jr.,
Resolutions with a given {H}ilbert function, in \emph{Commutative
algebra: syzygies, multiplicities, and birational algebra (South
Hadley, MA, 1992)}, in \emph{Contemp. Math.}  \textbf{159},
pp.~19--26 (Amer. Math. Soc., Providence, RI, 1994).

\bibitem[F]{Francisco} C. Francisco, Minimal graded Betti
numbers and stable ideals, \emph{Comm. Algebra}, \textbf{31} (2003),
no.~10, 4971--4987.

\bibitem[G-L-R-R]{GelvinLaVictoreReedRichert} M. Gelvin; P. LaVictore;
J. Reed; B. Richert, Betti numbers and simplicial complexes,
preprint.

\bibitem[G-H-S]{GeramitaHarimaShin} A.V. Geramita, T. Harima, and
Y.S. Shin, Extremal point sets and Gorenstein ideals, \emph{Adv.
Math.}, \textbf{152} (2000), no.~1, 78--119.

\bibitem[G-S]{M2} D.~R. Grayson and M.~E. Stillman, \emph{Macaulay 2,
a software system for research in algebraic geometry}.
\verb|http://www.math.uiuc.edu/Macaulay2/|.


\bibitem[H-R-W]{HerzogReinerWelker} J.  Herzog, V.  Reiner, and V.
  Welker, Componentwise linear ideals and Golod rings, \emph{Michigan
    Math. J.}, \textbf{46} (1999), no.~2, 211--223.

\bibitem[H]{Hulett} H. A. Hulett, Maximal Betti Numbers with a Given
Hilbert Function, \emph{Comm. in Alg.}, \textbf{21} (1993), no.~7,
2335--2350.

\bibitem[Ho]{Hochster} M. Hochster, Cohen-Macaulay rings,
  combinatorics, and simplicial complexes, in \emph{Ring theory, II
    (Proc.  Second Conf., Univ. Oklahoma, Norman, Okla., 1975)},
  171--223, in \emph{Lecture Notes in Pure and Appl. Math.},
  \textbf{26}, Dekker, New York, 1977.

\bibitem[M-N]{MiglioreNagel} J. Migliore; U. Nagel, Reduced
arithmetically Gorenstein schemes and simplicial polytopes with
maximal Betti numbers, \emph{Adv. Math.} \textbf{180} (2003), no.~1,
1--63.

\bibitem[P]{Pardue} K. Pardue, Deformation classes of graded modules
    and maximal Betti numbers, \emph{Illinois J.  of Math.},
    \textbf{40} (1996), 564--585.

\bibitem[Ro]{Rodriguez} M. Rodriguez, Ideals that attain a given
Hilbert function, \emph{Illinois J. Math.}, \textbf{44} (2000),
no.~4, 821--827.

\bibitem[R]{Richert} B.~P. Richert, Smallest graded {B}etti
numbers, \emph{J. Algebra} \textbf{244} (2001), no.~1, 236--259.

\bibitem[S]{Stanley}
R. Stanley,
\newblock \underline{Combinatorics and Commutative Algebra},
\newblock Burkh\"auser, Boston, 1983.

\end{thebibliography}
\end{document}